\documentclass[11pt]{amsart}

\usepackage{graphicx}

\newtheorem{theorem}{Theorem}[section]
\newtheorem{proposition}[theorem]{Proposition}
\newtheorem{lemma}[theorem]{Lemma}
\newtheorem{corollary}[theorem]{Corollary}

\newtheorem*{MinimumBridgeNumberTheorem}{Theorem~\ref{thm:minbridgenum}}
\newtheorem*{MaximumBridgeNumberTheorem}{Theorem~\ref{thm:max}}
\newtheorem*{BridgeNumberSetTheorem}{Theorem~\ref{thm:bridgenumberset}}
\newtheorem*{maxbridge}{Corollary~\ref{coro:max}}

\theoremstyle{definition}

\newtheorem{remark}[theorem]{Remark}

\numberwithin{equation}{section}

\newcommand{\bc}{\operatorname{bc}}
\newcommand{\br}{\operatorname{br}}

\newcommand{\depth}{\operatorname{depth}}
\newcommand{\dist}{\operatorname{dist}}
\newcommand{\length}{\operatorname{length}}

\newcommand{\Nbd}{\operatorname{Nbd}}

\newcommand{\rbc}{\operatorname{rbc}}

\renewcommand{\P}{\operatorname{{\mathcal P}}}

\newcommand{\D}{\operatorname{{\mathcal D}}}
\newcommand{\G}{\operatorname{{\mathcal G}}}
\newcommand{\K}{\operatorname{{\mathcal K}}}
\newcommand{\T}{\operatorname{{\mathcal T}}}

\begin{document}

\title[Tunnel leveling, depth, and bridge numbers]
{Tunnel leveling, depth, and bridge numbers}

\author{Sangbum Cho}
\address{Department of Mathematics\\
University of California at Riverside\\
Riverside, CA 92521\\
USA}
\email{scho@math.ucr.edu}

\author{Darryl McCullough}
\address{Department of Mathematics\\
University of Oklahoma\\
Norman, Oklahoma 73019\\
USA}
\email{dmccullough@math.ou.edu}
\urladdr{www.math.ou.edu/$_{\widetilde{\phantom{n}}}$dmccullough/}
\thanks{The second author was supported in part by NSF grant
DMS-0802424}


\date{\today}


\begin{abstract}
We use the theory of tunnel number $1$ knots introduced in \cite{CM} to
strengthen the Tunnel Leveling Theorem of Goda, Scharlemann, and Thompson.
This yields considerable information about bridge numbers of tunnel number
$1$ knots. In particular, we calculate the minimum bridge number of a knot
as a function of the maximum depth invariant $d$ of its tunnels.  The
growth of this value is on the order of $(1+\sqrt{2})^d$, which improves
known estimates of the rate of growth of bridge number as a function of the
Hempel distance of the associated Heegaard splitting. We also find the
maximum bridge number as a function of the number of cabling constructions
needed to produce the tunnel, showing in particular that the maximum bridge
number of a knot produced by $n$ cabling constructions is the $(n+2)^{nd}$
Fibonacci number. Finally, we examine the special case of the ``middle''
tunnels of torus knots.
\end{abstract}

\maketitle

\section*{Introduction}
\label{sec:intro}

The Tunnel Leveling Theorem of H. Goda, M. Scharlemann, and
A. Thompson~\cite{GST} says that when a tunnel number one knot is in
minimal bridge position, any of its tunnel arcs can be slid to lie in a
single horizontal level. Using the theory of tunnel number $1$ knots
developed in~\cite{CM}, we will prove the \textit{Tunnel Leveling
Addendum.}  Roughly speaking, it says that when a tunnel arc is in level
position as in the conclusion of the Tunnel Leveling Theorem, the other two
knots from the $\theta$-curve which is the union of the knot and its tunnel
arc are also (after trivial repositioning) in minimal bridge position. Its
full statement is given near the start of Section~\ref{sec:efficient}.

The Tunnel Leveling Addendum gives a great deal of information about bridge
numbers of tunnel number $1$ knots. Some of these applications involve the
\textit{depth} invariant, which is defined using the theory
from~\cite{CM}. The depth of a tunnel, $\depth(\tau)$, is somewhat similar
to the (Hempel) distance $\dist(\tau)$ (see J. Johnson
\cite{JohnsonBridgeNumber} and Y. Minsky, Y. Moriah, and S. Schleimer
\cite{MMS}), but unlike the distance, the depth is very easy to calculate
in terms of the parameter description of tunnels given in~\cite{CM}. The
two invariants are related by the inequality
\[\dist(\tau)-1 \leq \depth(\tau)\]
but the depth can be much larger than the distance. Indeed, the ``middle''
tunnels of torus knots that we examine below are easily seen to have
distance~$2$, but we will see their depths can be arbitrarily large.

The depth invariant can be defined geometrically in terms of the cabling
constructions of~\cite{CM}, but it also has a geometric interpretation in
terms of a construction that first appeared in~\cite{GST}. That
construction, which we call a giant step, is studied
in~\cite{CMgiantsteps}.

There is no upper bound for the bridge number of a knot in terms of the
depths of its tunnels, but among our applications of the Tunnel Leveling
Addendum is a sharp lower bound:
\begin{MinimumBridgeNumberTheorem}[Minimum Bridge Number]
For $d\geq 1$, the minimum bridge number of a knot having a tunnel of
depth~$d$ is given recursively by $a_d$, where $a_1=2$, $a_2=4$, and
$a_d=2a_{d-1}+a_{d-2}$ for $d\geq 3$. Explicitly,
\[a_d = \frac{(1+\sqrt{2})^d}{\sqrt{2}}- \frac{(1-\sqrt{2})^d}{\sqrt{2}}\]
and consequently ${\displaystyle \lim_{d\to\infty}} a_d -
\frac{(1+\sqrt{2})^d}{\sqrt{2}} = 0$.
\end{MinimumBridgeNumberTheorem}
This improves Lemma~2 of \cite{JohnsonBridgeNumber}, which is that bridge
number grows at least linearly with distance. It also improves
Proposition~1.11 of \cite{GST}, which implies that bridge number grows
at least as fast as~$2^d$.

Actually, the Minimum Bridge Number Theorem can be proven using only the
Tunnel Leveling Theorem, which give the lower bounds, and some explicit
constructions to realize the minimum values. Our upper bound result,
however, uses the full strength of the Tunnel Leveling Addendum:
\begin{MaximumBridgeNumberTheorem}[Maximum Bridge Number]
Write the Fibonacci sequence $(1,1,2,3,\ldots)$ as $(F_1,F_2,\ldots)$. The
maximum bridge number of a knot having a tunnel produced by $n$ cabling
constructions, of which the first $m$ produce simple or semisimple tunnels,
is $mF_{n-m+2}+F_{n-m+1}$.
\end{MaximumBridgeNumberTheorem}
\noindent 
(The terms ``simple'' and ``semisimple'' are recalled in
Section~\ref{sec:cabling}.) For fixed $n$, the largest value for the upper
bound in Theorem~\ref{thm:max} occurs when $m=2$, giving the following
absolute maximum:
\begin{maxbridge} 
The maximum bridge number of a knot
having a tunnel produced by $n$ cabling operations is~$F_{n+2}$.
\end{maxbridge}
\noindent
In fact, this maximum bridge number is achieved by a sequence of torus knot
tunnels, as we will see in Proposition~\ref{prop:depth-efficienttorusknots}.

In addition to giving general bounds, the Tunnel Leveling Addendum places
very strong restrictions on the possible bridge numbers that can occur:
\begin{BridgeNumberSetTheorem}[Bridge Number Set]
Suppose that a knot $K$ has a tunnel
$\tau$ produced by $n\geq 2$ cabling operations,
of which the first $m$ produce simple or semisimple tunnels. Then
$\br(K)$ is one of the $2m-2$ values $F_\tau(a,b)$ for $2\leq a\leq
b\leq a+1\leq m+1$.
\end{BridgeNumberSetTheorem}
\noindent Here, $F_\tau$ is the \textit{Fibonacci function} of $\tau$,
defined in Section~\ref{sec:Fibonacci_functions}. It appears almost certain
that all $2m-2$ possible values in the Bridge Number Set Theorem do occur
as bridge numbers. As explained in Remark~\ref{rem:bridge_number_set}, this
is easy to see for $m=2$ and $m=3$, but for the general case we have not
been able to verify all the necessary examples.

We will also examine the interesting case of the ``middle'' tunnels of
torus knots. In our paper~\cite{CMtorus}, we calculated the invariants
of~\cite{CM} for all torus knot tunnels. Using that information, we will
show that torus knot tunnels achieve the minimum rate of growth of bridge
numbers in terms of depth, but not the minimum possible values, while they
do achieve the maximum possible bridge numbers in terms of the number of
cabling constructions.

Here is an outline of the sections of the paper. The first two sections
constitute a concise review of the material from~\cite{CM} that we will
need for the present applications. Section~\ref{sec:ddd} introduces the
distance and depth invariants, and gives a few results that follow quickly
from~\cite{CM} and work of other authors. Section~\ref{sec:tunnel_leveling}
reviews the Tunnel Leveling Theorem, and Section~\ref{sec:efficient} states
and proves the Tunnel Leveling Addendum. Fibonacci functions are introduced
in Section~\ref{sec:Fibonacci_functions}, which also contains the more
technical results on bridge number, including the Bridge Number Set
Theorem. The Minimum and Maximum Bridge Number Theorems are proved in
Section~\ref{sec:growth}, and torus knot tunnels are studied in
Section~\ref{sec:torus_knots}. Finally, most of the results apply to the
case of tunnel number $1$ links, and in Section~\ref{sec:links},
we briefly discuss these adaptations.

\section{The disk complex of the genus-$2$ handlebody}
\label{sec:disk_complex}

\begin{figure}
\begin{center}
\includegraphics[width=4.5cm]{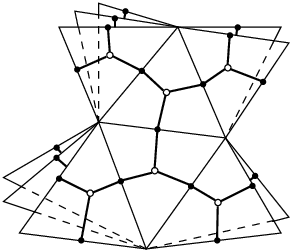}
\caption{A portion of the nonseparating disk complex $\D(H)$ and the tree
$\widetilde{\T}$. Countably many $2$-simplices meet along each edge.}
\label{fig:subdivision}
\end{center}
\end{figure}

Let $H$ be a genus~$2$ orientable handlebody, regarded as the standard
unknotted handlebody in $S^3$. For us, a \textit{disk in H} means a
properly imbedded disk in $H$, \textit{which is assumed to be nonseparating
unless otherwise stated.} The \textit{disk complex} $\D(H)$ is the
$2$-dimensional, contractible simplicial complex whose vertices are the
isotopy classes of disks in $H$, such that a collection of $k+1$ vertices
spans a $k$-simplex if and only if they admit a set of pairwise-disjoint
representatives. Each $1$-simplex of $\D(H)$ is a face of countably many
$2$-simplices. As suggested by Figure~\ref{fig:subdivision}, $\D(H)$ grows
outward from any of its $2$-simplices in a treelike way. In fact, it
deformation retracts to the tree $\widetilde{\T}$ seen in
Figure~\ref{fig:subdivision}.

Each disk $\tau$ in $H$ is the cocore disk of a tunnel of the knot $K_\tau$
which is a core circle of the solid torus obtained by cutting $H$ along
$\tau$. On the other hand, each tunnel of a tunnel number 1 knot in $S^3$
determines a collection of disks in $H$ as follows. The tunnel is a
$1$-handle attached to a regular neighborhood of the knot to form an
unknotted genus-$2$ handlebody. An isotopy carrying this handlebody to $H$
carries a cocore $2$-disk of that $1$-handle to a nonseparating disk in
$H$, and carries the tunnel number~$1$ knot to a core circle of the solid
torus obtained by cutting $H$ along the image disk in $H$. The
indeterminacy of this disk due to the choice of isotopy is the group of
isotopy classes of orientation-preserving homeomorphisms of $S^3$ that
preserve $H$. This group is called the \textit{Goeritz group $\G$.}  Work
of M. Scharlemann \cite{ScharlemannTree} and E. Akbas \cite{Akbas} proves
that $\G$ is finitely presented, and even provides a simple presentation of
it.

\begin{figure}\begin{center}
\includegraphics[width=4cm]{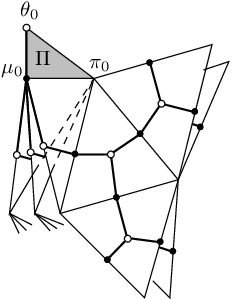}
\caption{A portion of $\D(H)/\G$ and $\T$ near the primitive orbits.}
\label{fig:Delta}
\end{center}
\end{figure}

Since two disks in $H$ determine equivalent tunnels exactly when they
differ by an isotopy moving $H$ through $S^3$, \textit{the collection of
all (equivalence classes of) tunnels of all tunnel number~$1$ knots
corresponds to the set of orbits of vertices of $\D(H)$ under $\G$.} So it
is natural to examine the quotient complex $\D(H)/\G$, which is illustrated
in Figure~\ref{fig:Delta}.
Through work of the first author~\cite{Cho}, the action of $\G$ on $\D(H)$
is well-understood. A disk $D$ in $H$ is called \textit{primitive} if there
is a disk $E$ in $\overline{S^3-H}$ for which $\partial D$ and $\partial E$
intersect transversely in one point in $\partial H$. The primitive disks
(regarded as vertices) span a contractible subcomplex $\P(H)$ of $\D(H)$,
called the \textit{primitive subcomplex}. The action of $\G$ on $\P(H)$ is
as transitive as possible, indeed the quotient $\P(H)/\G$ is a single
$2$-simplex $\Pi$ which is the image of any $2$-simplex of the first
barycentric subdivision of $\P(H)$. The vertices of $\Pi$ are $\pi_0$, the
orbit of all primitive disks, $\mu_0$, the orbit of all pairs of disjoint
primitive disks, and $\theta_0$, the orbit of all triples of disjoint
primitive disks.

On the remainder of $\D(H)$, the stabilizers of the action are as small as
possible. A $2$-simplex which has two primitive vertices and one
nonprimitive is identified with some other such simplices, then folded in
half and attached to $\Pi$ along the edge $\langle \mu_0,\pi_0\rangle$.
The nonprimitive vertices of such $2$-simplices are exactly the disks in
$\D(H)$ that are disjoint from some primitive pair, and these are called
$\textit{simple}$ disks. As tunnels, they are the upper and lower tunnels
of $2$-bridge knots. The remaining $2$-simplices of $\D(H)$ receive no
self-identifications, and descend to portions of $\D(H)/\G$ that are
treelike and are attached to one of the edges $\langle \pi_0,\tau_0\rangle$
where $\tau_0$ is simple.

The tree $\widetilde{\T}$ shown in Figure~\ref{fig:subdivision} is
constructed as follows.  Let $\D'(H)$ be the first barycentric subdivision
of $\D(H)$. Denote by $\widetilde{\T}$ the subcomplex of $\D'(H)$ obtained
by removing the open stars of the vertices of $\D(H)$. It is a bipartite
graph, with ``white'' vertices of valence $3$ represented by triples and
``black'' vertices of (countably) infinite valence represented by
pairs. The valences reflect the fact that moving along an edge from a
triple to a pair corresponds to removing one of its three disks, while
moving from a pair to a triple corresponds to adding one of infinitely many
possible third disks to a pair.  The possible disjoint third disks that can
be added are called the \textit{slope disks} for the pair.

The image $\widetilde{\T}/\G$ of $\widetilde{\T}$ in $\D'(H)/\G$ is a tree
$\T$. The vertices of $\D'(H)/\G$ that are images of vertices of $\D(H)$
are not in $\T$, but their links in $\D'(H)/\G$ are subcomplexes of
$\T$. These links are known to be infinite trees. For each such vertex
$\tau$ of $\D'(H)/\G$, i.~e.~each tunnel, there is a unique shortest path
in $\T$ from $\theta_0$ to \textit{the vertex in the link of $\tau$ that is
closest to $\theta_0$.} This path is called the \textit{principal path} of
$\tau$, and this closest vertex is a triple, called the \textit{principal
vertex} of $\tau$. The two disks in the principal vertex, other than
$\tau$, are called the \textit{principal pair} of $\tau$. They are exactly
the disks called $\mu^+$ and $\mu^-$ that play a key role
in~\cite{Scharlemann-Thompson}. Figure~\ref{fig:ppath} below shows the
principal path of a certain tunnel.

The white vertices of $\T$ correspond to unknotted $\theta$-curves in
$S^3$, up to isotopy. For a white vertex gives a triple of nonseparating
disks, dual to a $\theta$-curve in $H$ in which each arc crosses one of the
disks and not the others. These are exactly the unknotted $\theta$-curves,
in that a regular neighborhood is isotopic to $H$ which is part of a
Heegaard splitting of $S^3$. Two such $\theta$-curves in $H$ are isotopic
in $S^3$ exactly when they are equivalent under the Goeritz group, so the
white vertices of $\T$ give the isotopy classification.

\section{The cabling construction and the binary invariants}
\label{sec:cabling}

\begin{figure}
\begin{center}
\includegraphics[width=\textwidth]{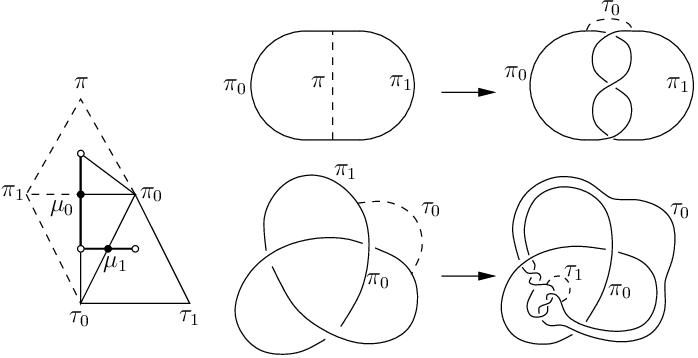}
\caption{Examples of the cabling construction.}
\label{fig:cabling}
\end{center}
\end{figure}

In a sentence, the cabling construction is to ``Think of the union of $K$
and the tunnel arc as a $\theta$-curve, and rationally tangle the ends of
the tunnel arc and one of the arcs of $K$ in a neighborhood of the other
arc of $K$.'' We sometimes call this ``swap and tangle,'' since one of the
arcs in the knot is exchanged for the tunnel arc, then the ends of other
arc of the knot and the tunnel arc are connected by a rational tangle.
Figure~\ref{fig:cabling} illustrates two cabling constructions, one
starting with the trivial knot and obtaining the trefoil, then another
starting with the tunnel of the trefoil.

More precisely, begin with a triple $\{\lambda,\rho,\tau\}$, where
$\{\lambda,\rho\}$ is the principal pair of $\tau$. Choose one of the disks
of $\{\lambda,\rho\}$, say $\rho$, and a slope disk $\tau'$ of the pair
$\{\rho,\tau\}$, \textit{other than $\lambda$.} This is a cabling operation
producing the tunnel $\tau'$ from $\tau$. The principal vertex of $\tau'$
is $\{\rho,\tau,\tau'\}$.

Unless otherwise stated, the slope disk $\tau'$ is chosen to be
nonseparating in $H$. A cabling operation using a separating disk as
$\tau'$ produces a tunnel number~$1$ link, and the cabling process cannot
be continued. This case will be discussed in Section~\ref{sec:links}.

Theorem~13.2 of~\cite{CM} shows that every tunnel of every tunnel
number~$1$ knot can be obtained by a uniquely determined sequence of
cabling constructions. A tunnel $\tau_0$ produced from the tunnel of the
trivial knot by a single cabling construction is called a \textit{simple}
tunnel. As already noted, these are the ``upper and lower'' tunnels of
$2$-bridge knots. A tunnel is called \textit{semisimple} if it is disjoint
from a primitive disk, but not from any primitive pair.

A $(1,1)$-knot is a knot that can be put into $1$-bridge position with
respect to a Heegaard torus of $S^3$. Let $K$ be a $(1,1)$-knot, whose
Heegaard torus splits $S^3$ into two solid tori $V$ and $W$. Associated to
this $(1,1)$-position are two tunnels obtained as follows. Let $\alpha_V$
be an arc in $V$ with endpoints in $K$, such that the union of $\alpha_V$
with the arc in $K\cap V$ bounded by the endpoints of $\alpha_V$ is a core
circle of $V$. Then $\alpha_V$ determines a tunnel of $K$; the
corresponding tunnel constructed in $W$ is the other one. Tunnels arising
in this way are called $(1,1)$-tunnels, and are exactly the simple and
semisimple tunnels.

A tunnel is called \textit{regular} if it is not primitive, simple, or
semisimple.

There is a procedure for assigning rational slopes which record the
rational tangle used in a cabling construction. We will not need these
slopes in our study of depth, although we will include them, in an
inessential way, in our discussion of the torus knot examples in
Section~\ref{sec:torus_knots}, and they also appear briefly in
Remark~\ref{rem:bridge_number_set}. The slope invariants are usually not
needed for working with depth because the depth is completely determined by
the second set of invariants associated to a tunnel, the ``binary''
invariants $s_2, s_3,\ldots\,$,~$s_n$, which we now define.

We have already mentioned that for every tunnel $\tau$, there is a unique
sequence of tunnels $\tau_0,\ldots\,$, $\tau_n=\tau$ such that $\tau_0$ is
simple and for each $i\geq 1$, $\tau_i$ is obtained from $\tau_{i-1}$ by a
cabling construction. The cabling that produces $\tau_i$ retains one arc of
the associated knot $K_{\tau_{i-1}}$ of $\tau_{i-1}$, and replaces the
other with a tangle, producing $K_{\tau_i}$. The invariant $s_i$ is $1$
exactly when this cabling replaces the arc that was retained by the
previous cabling, otherwise $s_i$ is $0$.

A tunnel is simple or semisimple if and only if all $s_i=0$. The reason is
that both conditions characterize cabling sequences in which one of the
original primitive disks is retained in every cabling; this corresponds to
the fact that the union of the tunnel arc and one of the arcs of the knot
is unknotted.

\begin{figure}
\begin{center}
\includegraphics[width=5cm]{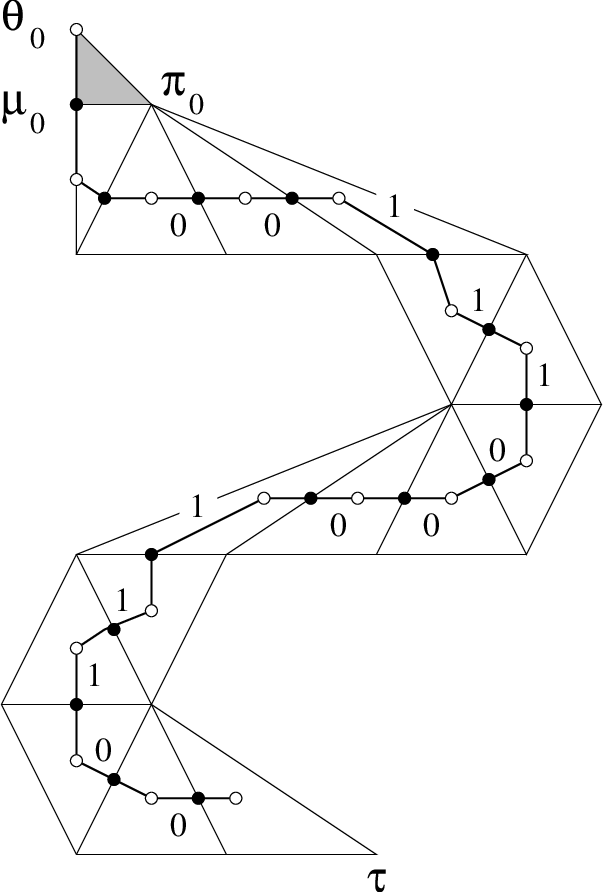}
\caption{The principal path of a tunnel $\tau$ having binary invariants
$0011100011100$, or equivalently with step sequence ``DRRRDRDLLLDLDRR''.}
\label{fig:ppath}
\end{center}
\end{figure}

There are two formal definitions of the binary invariants. The first is in
terms of the principal path $\theta_0$, $\mu_0$, $\mu_0\cup \{\tau_0\}$,
$\mu_1,\ldots\,$, $\mu_n$, $\mu_n\cup \{\tau_n\}$, where the $\mu_i$ are
the ``black'' vertices, the $\mu_i\cup\{\tau_i\}$ are the ``white''
vertices, and $\tau=\tau_n$: $s_i=0$ or $s_i=1$ according to whether or not
the unique disk of $\mu_i\cap\mu_{i-1}$ equals the unique disk of
$\mu_{i-1}\cap\mu_{i-2}$. Equivalently, each cabling operation begins with
a triple of disks $\{\lambda_{i-1},\rho_{i-1},\tau_{i-1}\}$ and finishes
with $\{\lambda_i,\rho_i,\tau_i\}$. For $i\geq 2$, put $s_i=1$ if
$\{\lambda_i,\rho_i,\tau_i\}=\{\tau_{i-2},\tau_{i-1},\tau_i\}$, and $s_i=0$
otherwise. Figure~\ref{fig:ppath} shows the principal path of a tunnel with
binary invariants $0011100011100$.

From the viewpoint of a traveler along the principal path, $s_i=1$ means a
change from making right turns (at the white vertices) to left turns, or
from left turns to right, while $s_i=0$ means a turn in the same direction
as the previous turn. Let us say that a \textit{step} of the principal path
is a portion between successive white vertices. A principal path can then
be described as a \textit{step sequence.} This is a string of symbols
``L'', ``R'', or ``D'', for ``left'', ``right'', and ``down'' as seen from
the reader's viewpoint (as opposed to the ``left'' and ``right'' of a
traveler along the path). For the example of Figure~\ref{fig:ppath}, the
step sequence is ``DRRRDRDLLLDLDRR''. In general, the initial step of a
principal path is always ``D'', and the second step, due to the standard
way that we draw the picture, is ``R''. Each subsequent step corresponds to
a binary invariant. An ``L'' can only be followed by another ``L'' or a
``D'', according as the corresponding binary invariant $s$ is $0$ or $1$,
and similarly an ``R'' is followed by another ``R'' or a ``D'', according
as $s$ is $0$ or $1$. When the previous step is ``D'', the effect of $s$
depends on the step before one that produces the ``D''. If the ``D'' is in
a sequence ``LD'', then the next step is ``R'' or ``L'' according as $s$ is
$0$ or $1$, while if it is in ``RD'', then the next step is ``L'' or ``R''
according as $s$ is $0$ or $1$.

Functions that translate between the binary sequence and step sequence
descriptions are included in the software at~\cite{slopes}. The main
functions there accept either form of input for principal paths.

\section{Distance and depth}
\label{sec:ddd}

The \textit{(Hempel) distance} $\dist(\tau)$ is the shortest distance in
the curve complex of $\partial H$ from $\partial \tau$ to a loop that
bounds a disk in $\overline{S^3-H}$ (see J. Johnson
\cite{JohnsonBridgeNumber} and Y. Minsky, Y. Moriah, and S. Schleimer
\cite{MMS}). It is well-defined since the action of the Goeritz group on
$\partial H$ preserves the set of loops that bound disks in $H$ and the set
that bound in~$\overline{S^3-H}$.

A nonseparating disk has distance $1$ if and only if it is primitive, since
both conditions are equivalent to the condition that cutting $H$ along the
disk produces an unknotted solid torus. Therefore the tunnel of the trivial
knot is the only tunnel of distance~$1$. A simple or semisimple tunnel has
distance~$2$, since it is disjoint from a primitive disk. There are,
however, regular tunnels of distance~$2$. It is an easy observation that
the ``middle'' tunnels of torus knots all have distance~$2$, and in most
cases these are regular.

When $\Sigma_\tau=(\overline{H-\Nbd(K_\tau)},\overline{S^3-H})$ is a
Heegaard splitting of the complement of $K_\tau$, the \textit{(Hempel)
distance $\dist(\Sigma_\tau)$} is the minimal distance in the curve complex
of $\partial H$ between the boundary of a disk in
$\overline{H-\Nbd(K_\tau)}$ and the boundary of a disk in
$\overline{S^3-H}$ (where the disks may be separating). Clearly,
$\dist(\Sigma_\tau)\leq \dist(\tau)$. On the other hand, Johnson
\cite[Lemma 11]{JohnsonBridgeNumber} proved that
\begin{lemma}[Johnson] $\dist(\tau)\leq \dist(\Sigma_\tau)+1$.
\label{lem:JohnsonDistance}
\end{lemma}

M. Scharlemann and M. Tomova~\cite{Scharlemann-Tomova} proved the following
stability result:
\begin{theorem}[Scharlemann-Tomova]
Genus-$g$ Heegaard splittings of distance more than $2g$ are isotopic.
\label{thm:Scharlemann-Tomova}
\end{theorem}
\noindent Using Lemma~\ref{lem:JohnsonDistance} and
Theorem~\ref{thm:Scharlemann-Tomova}, Johnson~\cite[Corollary
13]{JohnsonBridgeNumber} deduced the following:
\begin{theorem}[Johnson] If $\tau$ is a tunnel of a tunnel number $1$ knot
$K_\tau$ and $\dist(\tau)>5$, then $\tau$ is the unique tunnel of $K_\tau$.
\label{thm:uniquetunnel}
\end{theorem}

Theorem~15.2 of~\cite{CM} determines all orientation-reversing
self-equivalences of tunnels:
\begin{theorem} Let $\tau$ be a tunnel of a tunnel number~$1$ knot or link.
Suppose that $\tau$ is equivalent to itself by an orientation-reversing
equivalence. Then $\tau$ is the tunnel of the trivial knot, the trivial
link, or the Hopf link.
\label{thm:HopfLink}
\end{theorem}
Combining Theorems~\ref{thm:uniquetunnel} and~\ref{thm:HopfLink} gives the
following:
\begin{corollary}
If $\tau$ is a tunnel of a tunnel number $1$ knot, and $\dist(\tau)>5$,
then $K_\tau$ is not amphichiral.
\label{coro:amphichiral}
\end{corollary}
\noindent For Theorem~\ref{thm:HopfLink} shows that an
orientation-reversing equivalence from $K_\tau$ to $K_\tau$ would produce a
second tunnel for~$K_\tau$.

Distance also has implications for hyperbolicity:
\begin{theorem} If $K_\tau$ is a torus knot or a satellite knot, then
$\dist(\tau)\leq 2$.
Consequently, if $\dist(\tau)\geq 3$, then $K_\tau$ is hyperbolic.\par
\label{thm:Morimoto-Sakuma}
\end{theorem}
\begin{proof} We have already mentioned the fact~\cite{CMtorus}
that the middle tunnels of torus knots have distance~$2$. The other tunnels
of torus knots are simple or semisimple, so also have
distance~$2$. K. Morimoto and M. Sakuma \cite{Morimoto-Sakuma} found all
tunnels of tunnel number~$1$ satellite knots, showing in particular that
they are semisimple.
\end{proof}

The \textit{depth} of $\tau$ is the simplicial distance $\depth(\tau)$ in
the $1$-skeleton of $\D(H)/\G$ from $\tau$ to the primitive vertex $\pi_0$.
From the definitions, $\tau$ is primitive if and only if $\depth(\tau)=0$,
is simple or semisimple if and only if $\depth(\tau)=1$, and is regular if
and only if $\depth(\tau)\geq 2$.

The inequality
\[\dist(\tau)-1\leq \depth(\tau)\]
mentioned in the introduction is immediate from the definitions. On the
other hand, we have already noted that the middle tunnels of torus knots
have distance $2$, but we will see in Section~\ref{sec:torus_knots} that
their depths can be arbitrarily large.

In terms of the step sequence describing the principal path of a tunnel,
the depth is simply the number of D's that appear. One can, of course,
determine the depth directly from the binary invariants. A maximal block of
$1$'s in the binary word $s_2\cdots s_n$ has the following effect: its
first, third, fifth, and so on terms will produce a downward step,
increasing the depth, while the other terms correspond to horizontal steps,
keeping the same depth. This gives the following simple algorithm to
compute $\depth(\tau)$ from the binary invariants of $\tau$:
\begin{enumerate}
\item Write the binary word $s_2s_3\cdots s_n$ as $O_1Z_1O_2Z_2\cdots
O_kZ_k$, where $O_i$ and $Z_i$ are respectively maximal blocks of ones and
zeros (thus $O_1$ and $Z_k$ may have length $0$, while all others have
positive length).
\item The depth of $\tau$ is $1+\displaystyle\sum_{j=1}^k \lceil
\length(O_i)/2\rceil$, where $\lceil
\length(O_i)/2\rceil$ denotes the least integer greater
than or equal to $\length(O_i)/2$.
\end{enumerate}

\section{Tunnel leveling}
\label{sec:tunnel_leveling}

Roughly speaking, the Tunnel Leveling Theorem of Goda, Scharlemann, and
Thompson says that a tunnel arc of a tunnel number one knot can be slid so
that it lies in a level sphere of some minimal bridge position of the
knot. Here is the rather technical version of the Tunnel Leveling Theorem
that we will need. Illustrations of conclusions (i) and (ii) of the theorem
appear in the first drawings of Figure~\ref{fig:cc} and
Figure~\ref{fig:cc_eyeglass_semisimple} respectively.
\begin{theorem}[Goda-Scharlemann-Thompson] Let $\{\lambda,\rho\}$ be the 
principal pair of a tunnel $\tau$, and let $\theta$ be the $\theta$-curve
associated to the principal vertex $\{\lambda,\rho, \tau\}$ of
$\tau$. Write $T$ for the arc dual to $\tau$, and $L$ and $R$ for the other
two arcs of $\theta$ that are dual to $\lambda$ and $\rho$, so that
$K_{\tau} = L \cup R$, $K_\lambda=R \cup T$, and $K_\rho=L \cup T$. Then
there is a minimal bridge position of $K_\tau$ for which either:
\begin{enumerate}
\item[(i)] $T$ is slid to an arc in a level sphere, and $T$ connects two
bridges of $K_\tau$. Moreover, $K_\tau\cup T$
is isotopic to the original~$\theta$. Or,
\item[(ii)] $T$ is slid to an eyeglass in a level sphere. The endpoints of
$T$ can be slid slightly apart, moving $T$ out of the level sphere,
producing $K_\tau\cup T$ isotopic to the original $\theta$, and showing
that one of $K_\lambda$ or $K_\rho$ is a trivial knot, and consequently
$\tau$ is simple or semisimple.
\end{enumerate}
and furthermore, in the $n$-strand trivial tangle above the level
sphere:
\begin{enumerate}
\item[(iii)] In case (i), the arcs are parallel to a collection of disjoint
arcs in the level sphere, which meet $T$ only in its endpoints.
\item[(iv)] In case (ii), the $n-1$ arcs not meeting $T$ are parallel to a
collection of disjoint arcs in the level sphere, each meeting the eyeglass
in a single point.
\end{enumerate}
\label{thm:GSTtheorem}
\end{theorem}
\begin{proof}
By Theorem~1.8 of \cite{GST}, we may move $K_\tau \cup T$, possibly using
slide moves of $T$ as well as isotopy, so that $K_\tau$ is in minimal
bridge position and $T$ either lies on a level sphere and connects two
bridges of $K_\tau$, or $T$ is slid to an ``eyeglass''. Since the leveling
process involves sliding the tunnel arc $T$, there is \textit{a priori} no
reason for the resulting $\theta$-curve to be isotopic to the original
$\theta$. But Corollary~3.4 and Theorem~3.5 (combined with Lemma~2.9) of
\cite{Scharlemann-Thompson} show that in~(i) and~(ii), the dual disks to
the other two arcs of the $\theta$-curve are the disks called $\mu^+$ and
$\mu^-$ there. By Lemma~14.1 of \cite{CM}, these disks are
the principal pair of $\tau$, that is, $\lambda$ and $\rho$. Therefore
the resulting $\theta$-curve is isotopic
to the original~$\theta$. Finally, the description of the trivial tangle
above the level sphere in~(iii) is from Theorem~6.1 of~\cite{GST}, and
in~(iv) from Corollary~6.2 of~\cite{GST}, which relies on~\cite{GOT}.
\end{proof}

A tunnel arc $T$ satisfying conclusion (i) of the Tunnel Leveling Theorem
is said to be in \textit{level arc position,} while for conclusion~(ii),
\textit{after sliding the endpoints apart to produce $\theta$,} it is in
\textit{eyeglass position.} If it is in one of these two positions, it may
be said to be in \textit{level position.}

A tunnel arc satisfying all the requirements of level position except that
the number of bridges of $K_\tau$ is not necessarily minimal is said to be
in \textit{weak level position.} The number of bridges is then called the
\textit{bridge count,} denoted $\bc(K_\tau)$ and dependent, of course, on
the choice of weak level position.

Suppose that $\tau$ is in weak level arc position. The endpoints of $\tau$
cut $K_\tau$ into two arcs, one dual to $\lambda$ and the other dual to
$\tau$. By a simple isotopy, we may assume that one end of the arc dual to
$\lambda$ leaves the endpoints of the tunnel arc in the upward direction,
and the other leaves in the downward direction. For if both leave in the
same direction, we can slide an endpoint of the tunnel arc over one of the
arches, achieving a level position for which the two ends leave in
different directions. We then call this an \textit{admissible} weak level
arc position.

When $\tau$ is in weak level arc position, each of the local maxima of
$K_\tau$ lies in exactly one of $K_\lambda$ or $K_\rho$. The numbers that
lie in each are called the \textit{relative bridge counts} of $K_\lambda$
and $K_\rho$ for the weak level arc position of $\tau$, and denoted by
$\rbc(K_\lambda)$ and $\rbc(K_\rho)$. Clearly
$\bc(K_\tau)=\rbc(K_\lambda)+\rbc(K_\rho)$. If $\tau$ is in weak eyeglass
position, with $K_\lambda$ a trivial knot, then we define
$\rbc(K_\lambda)=1$ and $\rbc(K_\rho)=\bc(K_\tau)$. One always has
$\br(K_\gamma)\leq \rbc(K_\gamma)$.

A first consequence of Theorem~\ref{thm:GSTtheorem} is the following.
\begin{lemma} Let $\tau$ be a tunnel of a nontrivial knot, and
let $\{\lambda,\rho\}$ be the principal pair of $\tau$. Then
$\br(K_\lambda)+\br(K_\rho)-1\leq \br(K_\tau)$. If $\tau$ is regular,
then $\br(K_\lambda)+\br(K_\rho)\leq \br(K_\tau)$.
\label{lem:bridge_number_inequality}
\end{lemma}
\begin{proof} Apply Theorem~\ref{thm:GSTtheorem} to $\tau$.
If the tunnel arc is in level arc position, which may be assumed to be
admissible, then we have $\br(K_\lambda)+\br(K_\tau)\leq \rbc(K_\lambda) +
\rbc(K_\tau)=\br(K_\tau)$. If the tunnel arc is in eyeglass position,
producing, say, $K_\lambda$ trivial, then we have
$\br(K_\lambda)+\br(K_\rho)\leq 1 + \rbc(K_\rho) = 1 + \br(K_\tau)$, giving
the inequality. If $\tau$ is regular, then the eyeglass configuration
cannot occur, giving the stronger inequality.
\end{proof}

\section{Efficient cabling and the Tunnel Leveling Addendum}
\label{sec:efficient}

In this section, we will prove the following theorem:
\begin{theorem}[Tunnel Leveling Addendum] Let $\tau$ be a tunnel with
principal vertex $\{\lambda,\rho,\tau\}$. If $\tau$ is not simple, choose
notation so that $\rho$ is the tunnel directly preceding $\tau$ in the
cabling sequence for $\tau$. Assume that $\tau$ is not the tunnel of
the trivial knot or a simple tunnel of a $(2n+1,2)$ torus knot. Then either
\begin{enumerate}
\item[(a)] All level positions of $\tau$ are level arc positions, and
$\br(K_\tau)=\br(K_\rho) + \br(K_\lambda)$, or
\item[(b)] All level positions of $\tau$ are eyeglass positions, $\tau$ is
semismiple, and $\br(K_\tau)=\br(K_\rho)$.
\end{enumerate}
\label{thm:addendum}
\end{theorem}

\noindent The exceptional case of the Tunnel Leveling Addendum is detailed
in the next theorem, which is simply a restatement of some results from the
work of K. Morimoto and M. Sakuma on tunnels of $2$-bridge knots
\cite{Morimoto-Sakuma}:
\begin{theorem} The trefoil knots have unique tunnels, which are simple and 
can be put into either level arc position or eyeglass position. For the
other $(2n+1,2)$ torus knots, there are two simple tunnels, of which one
can be put into both kinds of level position, and the other only into level
arc position.
\label{thm:exceptional_case}
\end{theorem}
\noindent Of course the trivial knot has a unique tunnel, which can only be
leveled in eyeglass position.

It is important to understand the geometric content of the Tunnel Leveling
Addendum from the viewpoint of the Tunnel Leveling
Theorem~\ref{thm:GSTtheorem}. Apart from the exceptional cases, the
Addendum says that when a tunnel $\tau$ is leveled, giving a positioning of
the $\theta$-curve associated to the principal vertex
$\{\lambda,\rho,\tau\}$ of $\tau$, then (possibly after trivial
repositioning) the copies of $K_\lambda$ and $K_\rho$ in that
$\theta$-curve are in minimal bridge position.

\begin{figure}
\begin{center}
\includegraphics[width=\textwidth]{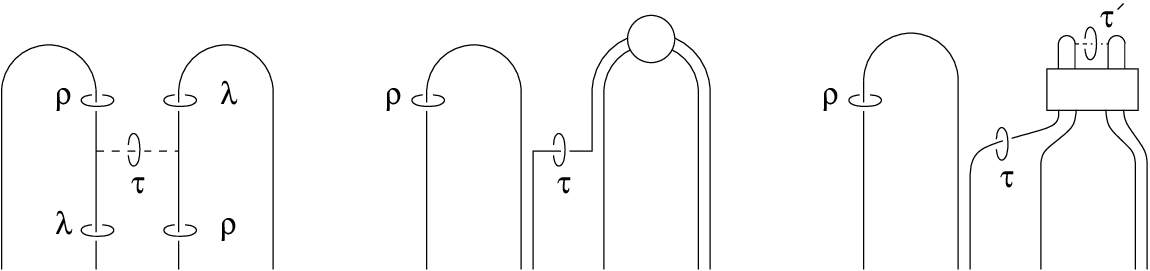}
\caption{Efficient cabling from admissible weak level arc position.}
\label{fig:cc}
\end{center}
\end{figure}

In this section we will prove the Tunnel Leveling Addendum and
Theorem~\ref{thm:exceptional_case}, and in preparation for this we now
introduce the technique of efficient cabling. The basic construction is
shown in Figure~\ref{fig:cc}, whre notation is selected so that the cabling
will replace $\lambda$ and retain $\rho$. We start with $\tau$ in
admissible weak level arc position, as shown in the left-hand drawing.
There may, of course, be many more bridges, some in $K_\rho$ and some
in~$K_\lambda$. A cabling of some arbitrary slope replaces $\lambda$ with
a new tunnel $\tau'$; the rational tangle in $K_{\tau'}$ created by the
cabling is inside a ball represented by the circle in the middle
drawing. We may then reposition $K_{\tau'}$ as in the right-hand drawing of
Figure~\ref{fig:cc}, by ``moving the ball up to engulf infinity,'' in such
a way that the rectangle in the drawing contains a $4$-strand braid. The
arc dual to $\tau'$ is in weak level arc position, and by a further
isotopy, if necessary, we may assume that it is in admissible weak level
arc position.

\begin{figure}
\begin{center}
\includegraphics[width=\textwidth]{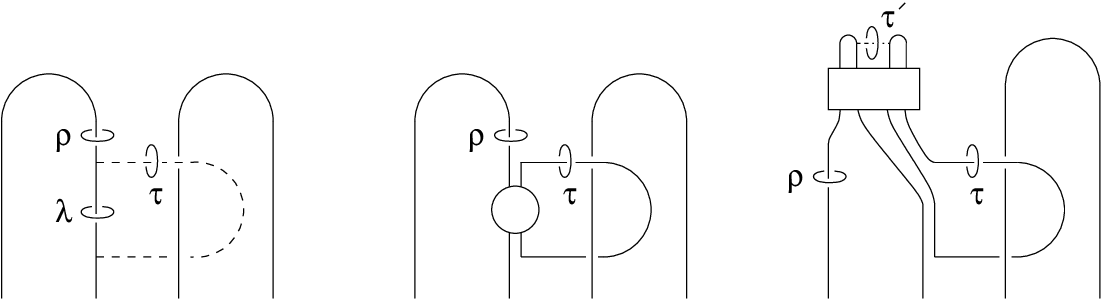}
\caption{Efficient cabling from eyeglass position, producing a
semisimple tunnel.}
\label{fig:cc_eyeglass_semisimple}
\end{center}
\end{figure}
\begin{figure}
\begin{center}
\includegraphics[width=\textwidth]{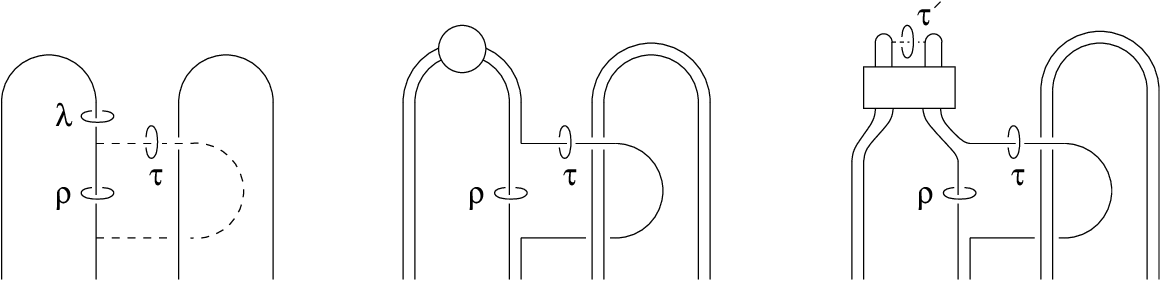}
\caption{Efficient cabling from eyeglass position, producing a
regular tunnel.}
\label{fig:cc_eyeglass_regular}
\end{center}
\end{figure}

The corresponding construction for a tunnel in weak eyeglass position can
produce either another semisimple tunnel or a regular tunnel. The resulting
tunnel is in weak level arc position, which by isotopy is also assumed to
be admissible. Efficient cablings for each of the two possibilities are
shown in Figures~\ref{fig:cc_eyeglass_semisimple}
and~\ref{fig:cc_eyeglass_regular}, and the constructions should be clear
from the discussion of the weak level arc case.

The next result details the effect of efficient cabling on bridge counts.
The notations $K$ and $K'$ are used to indicate knots obtained from the
$\theta$-curves associated to $\{\lambda,\rho,\tau\}$ and
$\{\rho,\tau,\tau'\}$ respectively:
\begin{proposition}
Suppose that $\tau$ is in weak level arc or weak eyeglass position, and 
that a cabling operation as in Figure~\ref{fig:cc},
\ref{fig:cc_eyeglass_semisimple}, or~\ref{fig:cc_eyeglass_regular} is
performed, producing a new tunnel $\tau'$ with principal vertex $\{\rho,
\tau, \tau'\}$, and producing a tunnel arc $T'$ for which $\tau'$ is in
level arc position. Then
\begin{enumerate}
\item $\rbc(K_\rho')=\rbc(K_\rho)$.
\item $\rbc(K_\tau')=\bc(K_\tau)$.
\item $\bc(K_{\tau'}') = \bc(K_\tau) + \rbc(K_\rho)$.
\end{enumerate}
\label{prop:bridge_counts}
\end{proposition}
\begin{proof}
The third equality follows from the first two. The first two are seen by
examination of Figures~\ref{fig:cc}, \ref{fig:cc_eyeglass_semisimple},
and~\ref{fig:cc_eyeglass_regular}. For example, let us consider
Figure~\ref{fig:cc}. In the leftmost drawing, denote the arcs dual to
$\lambda$, $\rho$, and $\tau$ by $L$, $R$, and $T$ respectively. In the
rightmost drawing, after the cabling producing $K_{\tau'}'$ has been
performed, denote the dual arcs by $R_1$, $T_1$, and $T_1'$, where the
latter is horizontal. By isotopy we may assume that $\tau'$ is also in
admissible weak level arc position.  The number of bridges that we see in
$R_1\cup T_1'$ equals the number that appeared in $R$ plus the number that
appeared in $L$, showing that $\rbc(K_\tau')=\bc(K_\tau)$. The number of
bridges in $T_1\cup T_1'$ is the number that appeared in $L$, so
$\rbc(K_\rho')=\rbc(K_\rho)$. The arguments for
Figures~\ref{fig:cc_eyeglass_semisimple} and~\ref{fig:cc_eyeglass_regular}
are similar.
\end{proof}

%\begin{proof}
%Examination of Figures~\ref{fig:cc}, \ref{fig:cc_eyeglass_semisimple},
%and~\ref{fig:cc_eyeglass_regular} verifies the first two equalities, and
%the third is an immediate consequence.
%\end{proof}

We are now ready to prove the Tunnel Leveling Addendum and
Theorem~\ref{thm:exceptional_case} simultaneously. The tunnel called $\tau$
in the statement of the Addendum will be denoted by $\tau'$ in our
argument. Its principal vertex will be written as $\{\rho,\tau,\tau'\}$. If
$\tau'$ is not simple, then we assume that $\tau$ is the tunnel that
precedes $\tau'$ in the cabling sequence of $\tau'$, and the principal
vertex of $\tau$ will be written as $\{\lambda,\rho,\tau\}$. As in
Proposition~\ref{prop:bridge_counts}, we use $K$ and $K'$ to indicate knots
obtained from the $\theta$-curves associated to $\{\lambda,\rho,\tau\}$ and
$\{\rho,\tau,\tau'\}$ respectively. Note, however, that $K_\rho$ and
$K_\tau$ are equivalent to $K'_\rho$ and~$K'_\tau$ and hence
$\br(K'_\rho)=\br(K_\rho)$ and $\br(K'_\tau)=\br(K_\tau)$.

We will induct on the length of the cabling sequence of $\tau'$. If the
length is $1$, then $\tau'$ is an upper or lower tunnel of the $2$-bridge
knot $K'_{\tau'}$, so can be put into level arc position.  Each of
$K'_\tau$ and $K'_\rho$ is a trivial knot so has bridge number
$1$. Therefore we have $\br(K'_{\tau'})=\br(K'_\rho) + \br(K'_\tau)$. So
conclusion~(a) of the Tunnel Leveling Addendum holds for $\tau'$, provided
that $\tau'$ cannot also be put into eyeglass position.

The homeomorphism classfication of tunnels of $2$-bridge knots is given in
Table 5.2(B) of \cite{Morimoto-Sakuma}, and only in the cases called $2$,
$3$, and $6$ there does there exist a simple tunnel which can also be put
into eyeglass position. Those cases are defined in Lemma~5.1 of
\cite{Morimoto-Sakuma}, and upon examination are found to be exactly the
$2$-bridge torus knots. Since this is the excluded case in the Tunnel
Leveling Addendum, the Addendum holds for tunnels whose cabling sequences
have length~$1$. Closer examination of the tunnel classification in
\cite{Morimoto-Sakuma} verifies the precise statement in
Theorem~\ref{thm:exceptional_case}, whose proof is now complete.

Assume now that the length of the cabling sequence of $\tau'$ is greater
than $1$. Put $\tau$ in level position and obtain $\tau'$ by efficient
cabling as in one of Figures~\ref{fig:cc},
\ref{fig:cc_eyeglass_semisimple}, or~\ref{fig:cc_eyeglass_regular}.

Suppose first that the resulting weak level arc position for $\tau'$ is
actually a level arc position. Using Proposition~\ref{prop:bridge_counts}
and induction, we have $\br(K'_{\tau'})=\rbc(K'_\tau)+\rbc(K'_\rho)=
\bc(K_\tau)+\rbc(K_\rho)=\br(K_\tau)+\br(K_\rho)$.

Suppose now that the weak level arc position for $\tau'$ is not level arc
position. Then $\br(K'_{\tau'})<\bc(K'_{\tau'})=
\rbc(K'_\tau)+\rbc(K'_\rho)=\bc(K_\tau)+\rbc(K_\rho)=
\br(K_\tau)+\br(K_\rho)$, so Lemma~\ref{lem:bridge_number_inequality} tells
us that $\tau'$ is semisimple and $\br(K_{\tau'}') =
\br(K_\tau)+\br(K_\rho)-1$. But $\rho$ is primitive, since the principal
vertex of every semismiple tunnel contains a primitive disk, and $\tau$ is
not primitive since $\tau'$ is not simple. Therefore $\br(K_\rho)=1$ and
$\br(K_\tau)= \br(K'_{\tau'})$.

In the latter case, $\tau'$ cannot be put into level arc position, since
then we would have $\br(K'_{\tau'})=\rbc(K'_\tau)+\rbc(K'_\rho)\geq
\br(K'_\tau)+\br(K'_\rho)> \br(K_\tau)$. So either all level positions are
level arc positions, or all are eyeglass positions.
This completes the proof of the Tunnel Leveling Addendum and
Theorem~\ref{thm:exceptional_case}.

The next two corollaries are convenient restatements of parts of the Tunnel
Leveling Addendum.
\begin{corollary}
Let $\tau$ be a regular tunnel with principal vertex
$\{\lambda,\rho,\tau\}$. Then $\br(K_\rho) + \br(K_\lambda) = \br(K_\tau)$.
\label{coro:reg_additivity}
\end{corollary}

\begin{corollary}
Let $\tau$ be a semisimple tunnel with principal vertex
$\{\pi_0,\rho,\tau\}$. Then $\br(K_\rho)= \br(K_\tau)$ or
$\br(K_\tau)=\br(K_\rho)+1$, according to whether all level positions for
$\tau$ are eyeglass positions or all are level arc positions.
\label{coro:ss_additivity}
\end{corollary}

\section{Fibonacci functions}
\label{sec:Fibonacci_functions}

Let $\tau$ be a tunnel and write the cabling sequence of $\tau$ as
$\tau_0$, $\tau_1\,\ldots$, $\tau_{m-1}$, $\tau_m,\ldots\,$,
$\tau_{n-1}=\tau$, where $\tau_i$ is simple or semisimple exactly when $i\leq
m-1$. That is, $\tau$ is produced by $n$ cablings, the first $m$ of which
produce depth-$1$ tunnels. In particular, $n=m$ when $\depth(\tau)=1$.

\begin{figure}
\begin{center}
\includegraphics[width=6cm]{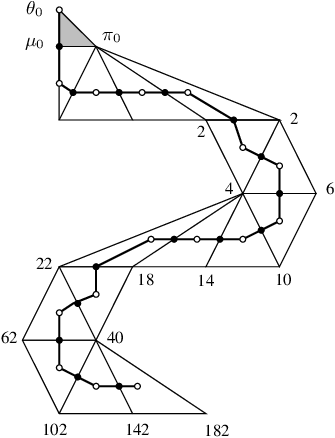}
\caption{}
\label{fig:bridge_numbers}
\end{center}
\end{figure}
The principal vertex of $\tau_0$ has the form $\{\pi_0,\pi_1,\tau_0\}$,
where $\pi_0$ and $\pi_1$ are primitive. If we put $\tau_{-1}=\pi_0$, the
trivial tunnel, then for each $k\geq 1$, the principal vertex of $\tau_k$
is of the form $\{\tau_i,\tau_{k-1},\tau_k\}$ for some $i$. If $n>m$, that
is, if $\tau$ is regular, then the first tunnel of depth $2$ is $\tau_m$,
and its principal vertex is $\{\tau_{m-2},\tau_{m-1},\tau_m\}$. We then
define the \textit{Fibonacci function} $F_\tau$ of a regular tunnel $\tau$
as follows. To compute $F_\tau(a,b)$, put $b_{m-2}=a$, $b_{m-1}=b$, and for
$m\leq k\leq n-1$, put $b_k=b_i+b_{k-1}$, where
$\{\tau_i,\tau_{k-1},\tau_k\}$ is the principal vertex of $\tau_k$.  Then,
put $F_\tau(a,b)=b_{n-1}$. Figure~\ref{fig:bridge_numbers} shows how to
calculate that $F_\tau(2,2)=182$ for a certain depth-$5$ tunnel with $m=4$
and $n=15$.

\begin{theorem}
Let $\tau$ be a simple or semisimple tunnel produced by $m$ cablings.
Then $2\leq \br(K_\tau) \leq m+1$.
\label{thm:ss_bridge_num}
\end{theorem}
\begin{proof}
Induct on $m$, using Corollary~\ref{coro:ss_additivity}.
\end{proof}

\begin{theorem}
Let $\tau$ be a regular tunnel whose cabling sequence contains $m$ tunnels
of depth $1$. Let $b_i=\br(K_{\tau_i})$ for $i\in \{m-2,m-1\}$. Then
$\br(K_\tau) =F_\tau(b_{m-2},b_{m-1})$.
\label{thm:reg_bridge_num}
\end{theorem}
\begin{proof}
Induct on the length of the cabling sequence of $\tau$, using
Corollary~\ref{coro:reg_additivity}.
\end{proof}

\begin{theorem}[Bridge Number Set]
Suppose that a knot $K$ has a tunnel
$\tau$ produced by $n\geq 2$ cabling operations,
of which the first $m$ produce simple or semisimple tunnels. Then
$\br(K)$ is one of the $2m-2$ values $F_\tau(a,b)$ for $2\leq a\leq
b\leq a+1\leq m+1$.
\label{thm:bridgenumberset}
\end{theorem}
\begin{proof}
By Theorem~\ref{thm:ss_bridge_num}, we have $2\leq \br(K_{\tau_{m-2}})\leq
m$, and by Corollary~\ref{coro:ss_additivity}, $\br(K_{\tau_{m-2}})\leq
\br(K_{\tau_{m-1}})\leq \br(K_{\tau_{m-2}})+1$. The result now follows from
Theorem~\ref{thm:reg_bridge_num}.
\end{proof}

\begin{remark}
We believe that for every principal path, each of the $2m-2$ possible
values given in the Bridge Number Set Theorem occurs as a bridge number for
some knots having a tunnel constructed using the given principal path. This
is clear for $m=2$. In that case, the cabling sequence has only two tunnels
$\tau_0$ and $\tau_1$ of depth $1$. When $\tau_1$ is a semisimple tunnel of
a $2$-bridge knot, $\br(K_{\tau_0})=\br(K_{\tau_1})=2$, and there are many
examples where $\br(K_{\tau_0})=2$ and $\br(K_{\tau_1})=3$, such as
semisimple tunnels of torus knots~\cite{CMtorus}. So choosing cabling
sequences that start with these two types of examples gives tunnels whose
knots have bridge numbers $F_\tau(2,2)$ and $F_\tau(2,3)$.

For $m=3$, we need to produce the sequences $(2,2,2)$, $(2,2,3)$,
$(2,3,3)$, and $(2,3,4)$ for
$(\br(K_{\tau_0}),\br(K_{\tau_1}),\br(K_{\tau_2}))$. Semisimple tunnels of
$2$-bridge knots give $(2,2,2)$. For $(2,2,3)$, we choose $\tau_1$ to be a
semisimple tunnel of a $2$-bridge knot and choose any cabling with slope
not of the form $\pm2 + 1/k$ to produce $\tau_2$; the results
of~\cite[Section 15]{CM} then show that $K_{\tau_2}$ cannot be
$2$-bridge. By Corollary~\ref{coro:ss_additivity}, $\br(K_{\tau_2})=3$. For
$(2,3,3)$, we start by constructing $\tau_1$ to be an upper tunnel of a
$3$-bridge torus knot, say the $(4,3)$ torus knot, as explained
in~\cite{CMtorus}. The tunnel arc shown in~\cite{CMtorus} can be put into
eyeglass level position with $K_{\tau_1}$ having three bridges, then a
cabling which is geometrically like those of Figure~14 of~\cite{CM} does
not raise bridge number, so $\br(K_{\tau_2})=3$ as well. Finally, for
$(2,3,4)$ we can just use the upper tunnel of the $(5,4)$ torus knot,
obtained by three cablings as in~\cite{CMtorus}.

For larger $m$, from upper tunnels of torus knots we obtain the bridge
number sequence $(2,3,4,\ldots,m+1)$, and hence realize $F_\tau(m,m+1)$.
And the idea for $(2,2,3)$ extends to realized $(2,2,\ldots,2,3)$, so
$F_\tau(2,3)$ occurs as a bridge number.  If we start with such an upper
tunnel sequence, for which the upper tunnel is in eyeglass position, and
then at some point begin using cablings as in Figure~14 of~\cite{CM}, we
obtain all sequences of the form $(2,3,4,\ldots,k-1,k,k,\ldots,k)$, giving
the $m-1$ values $F_\tau(k,k)$. So at least these $m+1$ values in the
bridge number set are known to occur. If we follow the latter procedure,
but use a complicated tangle for the final cabling, then the sequences
$(2,3,4,\ldots,k-1,k,k,\ldots,k,k+1)$ should also be obtained, giving the
remaining $m-3$ values $F_\tau(k,k+1)$ for $3\leq k\leq m-1$. Unfortunately
we lack a means to prove that the final knot has bridge number $k+1$ rather
than~$k$.
\label{rem:bridge_number_set}
\end{remark}

A peculiar consequence of Theorem~\ref{thm:reg_bridge_num} is the
following:
\begin{corollary}
Let $\tau$ be a regular tunnel and let $\tau_m$ be the first tunnel of
depth $2$ in the cabling sequence of $\tau$.  Then $\br(K_\tau)$ is
completely determined by the principal path of $\tau$ and the value of
$\br(K_{\tau_m})$. In fact,
\[\br(K_\tau)=F_\tau([\br(K_{\tau_m})/2], [(\br(K_{\tau_m})+1)/2] )\ .\]
\par
\label{coro:first_depth_2}
\end{corollary}
\begin{proof}
By Corollary~\ref{coro:reg_additivity}, $\br(K_{\tau_{m-2}}) +
\br(K_{\tau_{m-1}}) = \br(K_{\tau_m})$. Since also $\br(K_{\tau_{m-2}})$ and
$\br(K_{\tau_{m-1}})$ differ by at most $1$, their values must be as in the
statement of the corollary.
\end{proof}

It is not difficult to implement Theorems~\ref{thm:reg_bridge_num}
and~\ref{thm:bridgenumberset} computationally~\cite{slopes}:\par
\medskip

\begin{ttfamily}

\noindent Depth> fibonacci( '0011100011100', 2, 2, verbose=True )

\noindent F\underbar{\hspace*{0.8ex}}$\backslash$tau( 2, 2 ) = 182

\noindent The iteration sequence is:

\noindent \ \ \ 2, 2, 4, 6, 10, 14, 18, 22, 40, 62, 102, 142, 182

\noindent Depth> bridgeSet( '0011100011100' )

\noindent [182, 232, 273, 323, 364, 414]

\end{ttfamily}
\medskip

\section{Bounding bridge number}
\label{sec:growth}

Using the results of Section~\ref{sec:Fibonacci_functions}, we can give some
general bounds on bridge number. 

\begin{figure}
\begin{center}
\includegraphics{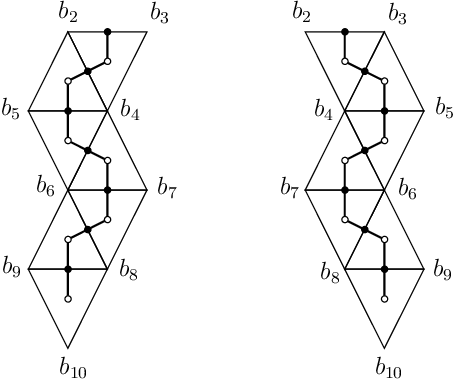}
\caption{The path on the left is of cheapest descent. The path on the right
is of cheapest descent if $b_2=b_3$.}
\label{fig:descent}
\end{center}
\end{figure}
First we examine lower bounds of bridge number as a function of depth. A
bit of experimentation with Fibonacci functions shows that the fastest
growth of depth relative to bridge number occurs for principal paths whose
regular portions (i.~e.~the parts starting from $\tau_{m-2}$ and
$\tau_{m-1}$) are the ``paths of cheapest descent'' seen in
Figure~\ref{fig:cheapest}. In that figure, the path on the left is always
cheapest descent, and the one on the right is cheapest descent when
$b_2=b_3$. Any principal path having more than two tunnels at a given depth
will produce a larger bridge number, as will any principal path that
emerges in the more costly direction out of a downward-pointing
$2$-simplex. From Theorem~\ref{thm:reg_bridge_num} we now have:
\begin{corollary} Let $\tau$ be a regular tunnel of depth $d$, and in
the principal path of $\tau$, let $\tau_m$ be the first tunnel of depth
$2$, with principal vertex $\{\tau_{m-2},\tau_{m-1},\tau_m\}$. Put
$b_2=\br(K_{\tau_{m-2}})$ and $b_3=\br(K_{\tau_{m-1}})$. For $n\geq 2$ let
$b_j$ be given by the recursion
\begin{gather*} b_{2n}=b_{2n-1}+b_{2n-2}\\
b_{2n+1}=b_{2n}+b_{2n-2}
\end{gather*}
Then $\br(K_\tau)\geq b_{2d}$.
\label{coro:cheapest_descent}
\end{corollary}

We can now prove one of our main results.
\begin{theorem}[Minimum Bridge Number]
For $d\geq 1$, the minimum bridge number of a knot having a tunnel of
depth~$d$ is given recursively by
$a_d$, where $a_1=2$, $a_2=4$, and $a_d=2a_{d-1}+a_{d-2}$ for
$d\geq 3$. Explicitly,
\[a_d = \frac{(1+\sqrt{2})^d}{\sqrt{2}}- \frac{(1-\sqrt{2})^d}{\sqrt{2}}\ .\]
and consequently ${\displaystyle \lim_{d\to\infty}} a_d -
\frac{(1+\sqrt{2})^d}{\sqrt{2}} = 0$.
\label{thm:minbridgenum}
\end{theorem}
\begin{proof} 
The smallest possible values for $\br(K_{\tau_{m-2}})$ and
$\br(K_{\tau_{m-1}})$ in Corollary~\ref{coro:cheapest_descent} are
$2$. These occur for any $m$, since there are $2$-bridge knot tunnels with
arbitrarily long cabling sequences, as seen in~\cite[Section~15]{CM}.
Taking $b_2=b_3=2$ in Corollary~\ref{coro:cheapest_descent} gives a
$b_{2d}$ which is a general lower bound for the bridge number of a tunnel
at depth~$d$, and a little bit of algebra shows that $b_{2d}=a_d$ for the
recursion in Theorem~\ref{thm:minbridgenum}. As a matrix, the recursion is
\[\begin{bmatrix} a_{d+1}\\
a_d\end{bmatrix} =
\begin{bmatrix} 2 & 1 \\
1 & 0
\end{bmatrix}
\begin{bmatrix} a_d\\
a_{d-1}
\end{bmatrix}\ .
\]
The eigenvalues of this matrix are $1\pm \sqrt{2}$, and elementary linear
algebra gives the formula $a_d = \frac{(1+\sqrt{2})^d}{\sqrt{2}}-
\frac{(1-\sqrt{2})^d}{\sqrt{2}}$.
\end{proof}

\begin{figure}
\begin{center}
\includegraphics[width=5cm]{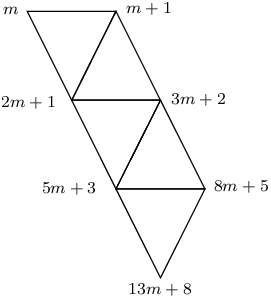}
\caption{The fastest growing upper bounds for bridge number, starting with
the last two semisimple tunnels in the cabling sequence.}
\label{fig:universal}
\end{center}
\end{figure}

We turn now to upper bounds. There is no upper bound in terms of depth,
since there are depth $1$ tunnels with arbitrarily large bridge number,
such as semisimple tunnels of torus knots~\cite{CMtorus}. We can, however,
bound the bridge number of $K_\tau$ in terms of the number of cablings
needed to produce~$\tau$. This time, we use the principal path forced by
choosing the larger of its two possible sums at every step, shown in
Figure~\ref{fig:universal}.
\begin{theorem}[Maximum Bridge Number]
Write the Fibonacci sequence $(1,1,2,3,\ldots)$ as $(F_1,F_2,\ldots)$. The
maximum bridge number of a knot having a tunnel produced by $n$ cabling
constructions, of which the first $m$ produce simple or semisimple tunnels,
is $mF_{n-m+2}+F_{n-m+1}$.
\label{thm:max}
\end{theorem}
\begin{proof}
If $\tau$ is simple, then $\br(K_\tau)=2$, $m=n=1$ and the expression
$mF_{n-m+2}+F_{n-m+1}$ equals $2$. If $\tau$ is semisimple, then $m=n$ and
$mF_{n-m+2}+F_{n-m+1}$ equals $m+1$, the upper bound given in
Theorem~\ref{thm:ss_bridge_num}. So we may assume that $\tau$ is regular.

In Figure~\ref{fig:universal}, the top two vertices are $\tau_{m-2}$ and
$\tau_{m-1}$, the last two semisimple tunnels that appear in the cabling
sequence of $\tau_{n-1}$. There are semisimple tunnels $\tau_{m-1}$
produced by $m$ cabling constructions which have $\br(K_{\tau_{m-1}})=m+1$,
such as the semisimple tunnels of the $(m+2,m+1)$ torus
knot~\cite{CMtorus}. Therefore the maximum bridge number is that given by
Theorem~\ref{thm:reg_bridge_num} applied to the principal path whose
regular portion is shown in Figure~\ref{fig:universal}. Using the fact
that $m =m\cdot F_1$ and $m+1=m\cdot F_2 + F_1$, one checks that this value
is $mF_{n-m+2}+F_{n-m+1}$.
\end{proof}

\begin{corollary} The maximum bridge number of a knot
having a tunnel produced by $n$ cabling operations
is~$F_{n+2}$.\par
\label{coro:max}
\end{corollary}
\begin{proof}
For a fixed $n$, the largest upper bound in Theorem~\ref{thm:max} occurs
when $m=2$.
\end{proof}

Proposition~\ref{prop:depth-efficienttorusknots} below gives an explicit
sequence of tunnels of torus knots that achieves the maximum value of
Corollary~\ref{coro:max}.

\section{Middle tunnels of torus knots}
\label{sec:torus_knots}

The tunnels of torus knots were classified by M. Boileau, M. Rost, and
H. Zieschang \cite{B-R-Z} and independently by Y. Moriah~\cite{Moriah}.

For a $(p,q)$ torus knot $K_{p,q}$ contained in the standard torus $T$ in
$S^3$, the \textit{middle tunnel} is represented by an arc in $T$ that
meets $K_{p,q}$ only in its endpoints. There are as many as two other
tunnels, which always have depth~$1$, but here we focus on the middle
tunnels.

In this section, we will include some information on slope invariants, for
those familiar with them. Slopes are not essential to the discussion, and
can be ignored if the reader so chooses.

For the tunnels of torus knots, the slope and binary invariants were
calculated in~\cite{CMtorus}. In particular, for the middle tunnels, we
have the following theorem, in which
$U=\begin{bmatrix}1&1\\0&1\end{bmatrix}$ and
$L=\begin{bmatrix}1&0\\1&1\end{bmatrix}$:
\begin{theorem} Let $p$ and $q$ be relatively prime integers with $p>q\geq
2$. Write $p/q$ as a continued fraction $[n_1,n_2,\ldots, n_k]$ with all
$n_j$ positive and $n_k\neq 1$. Let
\[ A_i = \begin{cases}%
L & -n_1\leq i\leq -1\\
U & 0\leq i\leq n_2-1\\
L & n_2\leq i\leq n_2+n_3-1\\
U & n_2+n_3\leq i\leq n_2+n_3+n_4-1\\
& \cdots \\
L& k\text{\ odd and\ }n_2+n_3+\cdots +n_{k-1}\leq i\leq n_2+n_3+\cdots
+ n_k-1\\
U& k\text{\ even and\ }n_2+n_3+\cdots +n_{k-1}\leq i\leq n_2+n_3+\cdots
+ n_k-1\ .
\end{cases}
\]
Put $N=n_2+n_3+\cdots + n_k-2$, and for $0\leq t\leq N$ put
\[\begin{bmatrix} a_t & b_t\\ c_t & d_t\end{bmatrix} =
\prod_{i=t}^{-n_1}A_i\ ,\]
where the subscripts in the product occur in descending order. Then:
\begin{enumerate}
\item[(i)] The middle tunnel of $K_{p/q}$ is
produced by $N+1$ cabling constructions whose slopes $m_0$, $m_1,\ldots\,$,
$m_N$ are
\[ \left[\frac{1}{2n_1+1}\right], \;a_1d_1+b_1c_1,\; a_2d_2+b_2c_2,\; \ldots,\;
a_Nd_N+b_Nc_N\ .\]
\item[(ii)] For each $t$, the cabling corresponding to the slope invariant
$m_t$ produces the $(a_t+c_t,b_t+d_t)$ torus knot; in particular, the first
cabling produces the $(2n_1+1,2)$ torus knot.
\item[(iii)] The binary invariants of the cabling sequence of this tunnel,
for $2\leq t\leq N$, are given by $s_t=1$ if $A_t\neq A_{t-1}$ and $s_t=0$
otherwise.
\end{enumerate}
\label{thm:middle_tunnels_slopes}
\end{theorem}
\noindent Note that this enables one to find the invariants of the middle
tunnels for all $(p,q)$ torus knots, since $K_{p,q}$ is isotopic to
$K_{q,p}$ and $K_{p,-q}$ is equivalent to $K_{p,q}$ by an
orientation-reversing homeomorphism taking middle tunnel to middle
tunnel. Such an equivalence negates the slope invariants and does not
change the binary invariants.

A bit of examination of the binary invariants yields a simple algorithm to
find the depth of the middle tunnel of $K_{p,q}$, $p> q\geq 2$:
\begin{enumerate}
\item Write $p/q$ as a continued fraction $[n_1,n_2,\ldots, n_k]$ with all
$n_i$ positive and $n_k\neq 1$.
\item Write the string $n_2\cdots n_k$ as $B_1B_2\cdots B_\ell$, where
each $B_i$ is either $n_in_{i+1}$ with $n_i=1$, or $n_i$ with $n_i\neq 1$.
\item The depth of the middle tunnel is $1+\ell$.
\end{enumerate}
This is implemented in the software at~\cite{slopes}.

\begin{figure}
\begin{center}
\includegraphics[height=9cm]{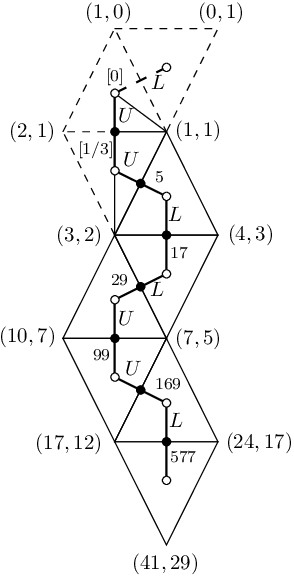}
\caption{Slowest growth of bridge number as a function of depth for torus
knot tunnels, corresponding to the continued fraction expansion
$41/29=[1,2,2,2,2]$. The $(41,29)$ torus knot has the smallest bridge
number of any torus knot with a depth~$4$ tunnel.}
\label{fig:cheapest}
\end{center}
\end{figure}
Figure~\ref{fig:cheapest} shows an initial segment of the principal paths
for the tunnels of the $(p,q)$-torus knots having continued fraction
expansions of the form $p/q=[1,2,2,\ldots,2]$. Notice that this is the path
of cheapest descent from Figure~\ref{fig:descent}. The small numbers along
the path are the slopes, the letters indicate whether the constructions
correspond to multiplication by $U$ or by $L$, and the pairs show the
$(p,q)$ for the torus knots determined by the tunnel at each step. The
first nontrivial cabling, with $m_0=[1/3]$, produces a $(3,2)$-torus knot,
and the second produces a $(4,3)$-torus knot with bridge number $3$. Since
we always have $p>q$, the bridge number is simply the value of $q$. These
obey the recursion of Corollary~\ref{coro:cheapest_descent}, starting with
$b_2=2$ and $b_3=3$. Since the cabling sequence for the middle tunnel
$\tau$ of any torus knot contains only one two-bridge knot (the
$(2n_1+1,2)$-torus knot produced by the first nontrivial cabling), there is
no regular torus knot tunnel which has $b_3=2$. Since the other tunnels of
torus knots are semisimple, the maximum depth of any tunnel of a torus knot
is the depth of its middle tunnel. Therefore each $b_{2d}$ in this sequence
gives the minimum bridge number for a torus knot with a tunnel of
depth~$d$. This gives a version of the Minimum Bridge Number
Theorem~\ref{thm:minbridgenum} for torus knot tunnels:
\begin{theorem}
For $d\geq 1$, the minimum bridge number of a torus knot tunnel
depth~$d$ is given recursively by
$t_d$, where $t_1=2$, $t_2=5$, and $t_d=2t_{d-1}+t_{d-2}$ for
$d\geq 3$. Explicitly,
%\[b_d = \frac{2+\sqrt{2}}{4}(1+\sqrt{2})^d +\
\[t_d = \frac{1}{2\sqrt{2}}(1+\sqrt{2})^{d+1} -
\frac{1}{2\sqrt{2}}(1-\sqrt{2})^{d+1}\ ,\]
and consequently ${\displaystyle \lim_{d\to\infty}} t_d -
\frac{1}{2\sqrt{2}}(1+\sqrt{2})^{d+1}=0$.
\end{theorem}
\noindent
We note that $t_d$ is $a_{d+1}/2$, where $a_{d+1}$ is the lower bound in
the Minimum Bridge Number Theorem~\ref{thm:minbridgenum}. That is, the
minimum bridge number of a torus knot having a tunnel of depth~$d$ is
exactly half the minimum bridge number for all knots having a tunnel of
depth $d+1$, and is approximately $(1+\sqrt{2})/2$ times the minimum
for all knots having a tunnel of depth~$d$.

In fact, the middle tunnel any torus knot for which $p/q$ has an expansion
$[n_1,2,2,2,2,\ldots,2]$ will have a principal path as in the previous
argument, since the first term in the continued fraction has no effect on
the principal path. Middle tunnels for which the expansion is not of this
form will have different principal paths, so we can state the following
result:
\begin{proposition}
The slowest growth of bridge number compared to depth for sequences of
middle tunnels of torus knots occurs when the $p/q$ have 
continued fraction expansions of the form $\pm[n_1,2,2,2,\ldots,2]$.
\label{prop:depth-efficienttorusknots}
\end{proposition}

\noindent By similar considerations, one can obtain the upper bound
version.
\begin{proposition}
The fastest growth of bridge number of torus knots per number of cablings
of the middle tunnels occurs for sequences of tunnels $\tau_k$ of $K_{p,q}$
for which the continued fraction expansions of $p/q$ are of the form
$\pm[n_1,1,1,1,\ldots,1]$, where there are $k$ $1$'s. For these tunnels,
$K_{\tau_k}$ has bridge number $F_{k+2}$.
\label{prop:fast-growthtorusknots}
\end{proposition}
\noindent For these tunnels, the terminal part of the corridor is like that
shown in Figure~\ref{fig:universal} with $m=2$. Since there are exactly $k$
cablings in the cabling sequence of $\tau_k$, these tunnels achieve the
bridge numbers in the Maximum Bridge Number Theorem. 

\section{The case of tunnel number $1$ links}
\label{sec:links}

As explained in ~\cite{CM}, our entire theory of tunnel number $1$ knot
tunnels can be adapted to include tunnels of tunnel number $1$ links, simply
by adding the separating disks as possible slope disks. The full disk
complex $\K(H)$ is only slightly more complicated than $\D(H)$. Each
separating disk is disjoint from only two other disks, both nonseparating,
so the additional vertices appear in $2$-simplices attached to $\D(H)$
along the edge opposite the vertex that is a separating disk. The quotient
$\K(H)/\G$ has only three types of additional $2$-simplices:
\begin{enumerate}
\item[(1)] There is a unique
orbit $\sigma_0$ of ``primitive'' separating disks, consisting of
separating disks disjoint from a primitive pair, which are exactly the
intersections of splitting spheres (see~\cite{ScharlemannTree}) with
$H$. In $\K(H)/\G$, $\sigma_0$ is a vertex of a ``half-simplex'' $\langle
\sigma_0, \pi_0, \mu_0\rangle$ attached to $\D(H)/\G$ along $\langle \pi_0,
\mu_0\rangle$. It is the unique tunnel of the trivial $2$-component link.
\item[(2)] Simple separating disks lie in
half-simplices attached along $\langle \pi_0, \mu_0\rangle$, just like
nonseparating simple disks.
\item[(3)] The remaining separating disks lie in $2$-simplices attached along
edges of $\D(H)/\G$ spanned by two (orbits of) disks, at least one of which
is nonprimitive. 
\end{enumerate}

For the spine, a single ``Y'' is added to $\T$ for each added $2$-simplex
as in~(3), and a folded ``Y'' for each the half-simplices as in~(1)
and~(2). The link in $\K'(H)/\G$ of a link tunnel is simply the top edges
(or top edge, for the trivial and simple tunnels) of such a~``Y'' (or
folded ``Y'').

The cabling operation differs only in allowing a separating slope disk, which
produces a tunnel of a tunnel number~$1$ link. The cabling sequence ends
with the first separating slope disk. Thus the principal paths look exactly
like those of the knot case, such as the one in Figure~\ref{fig:ppath}. The
only difference is that no further continuation is possible if the final
tunnel is the tunnel of a link.

For link tunnels, the distance and depth invariants are defined as for knot
tunnels. Simple tunnels are the upper and lower tunnels of $2$-bridge links
(and are the only tunnels of these links, see \cite{Adams-Reid},
\cite{Kuhn}, or~\cite[Theorem 16.3]{CM}). Depth~$1$ tunnels are the tunnels
of links with one component unknotted. The other component must be a
$(1,1)$-knot, and the link must have torus bridge number~$2$~\cite[Theorem
16.4]{CM}. Lemma~\ref{lem:JohnsonDistance} holds when $\tau$ is
separating, in fact the argument is an easier version of the argument
in~\cite{JohnsonBridgeNumber}, so Theorem~\ref{thm:uniquetunnel} and
Corollary~\ref{coro:amphichiral} hold for links as well as knots.

The Tunnel Leveling Addendum extends to tunnels of tunnel number $1$ links,
since the efficient cabling construction of Section~\ref{sec:efficient}
works just as well in the link case. But the statement and proof are very
much simpler, since only the level arc case need be considered.
\begin{theorem}[Tunnel Leveling Addendum for Links] Let $\tau$ be a 
tunnel of a tunnel number $1$ link, with principal vertex
$\{\lambda,\rho,\tau\}$. Then $\br(K_\tau)=\br(K_\rho) + \br(K_\lambda)$.
\label{thm:addendum_links}
\end{theorem}

Using the Tunnel Leveling Addendum for Links, the computational results of
Sections~\ref{sec:Fibonacci_functions} and~\ref{sec:growth} hold as stated
for tunnels of tunnel number $1$ links. Consequently, the software
implementations of \cite{slopes} also produce correct results for link
tunnels.

\bibliographystyle{amsplain}

\end{document}